\documentclass[12pt]{amsart}
\usepackage[centertags]{amsmath}

\usepackage{amsfonts}
\usepackage{amssymb}
\title[The pentagons and PA]{The pentagon as a substructure lattice of models of  Peano Arithmetic}
\author{James H. Schmerl}
\date{December~28,~2022}

\def\ta{{\sf TA}}

\def\pa{{\sf PA}}
\def\pas{\pa^*}
\def\PA{{\sf PA}}
\def\aca{{\sf ACA}_0}
\def\into{\longrightarrow}

\def\bbone {1 {\hspace{-4.6pt} 1}}
\def\bbzero{0 {\hspace{-4.6pt} 0}}

\DeclareMathOperator{\Lt}{Lt}
\DeclareMathOperator{\Ltr}{Ltr}
\DeclareMathOperator{\Eq}{Eq}
\DeclareMathOperator{\eq}{Eq}
\DeclareMathOperator{\Th}{Th}

\DeclareMathOperator{\Def}{Def}

\DeclareMathOperator{\ssy}{SSy}

\def\.{.\@}

\def\trieq{\unlhd}

\def\cal{\mathcal}

\def\MM{{\cal M}}
\def\NN{{\cal N}}
\def\KK{{\cal K}}

\def\res{\hspace{-4pt} \upharpoonright \hspace{-3pt}}
\def\harp{\res}

\def\id{{\sf id}}

\begin{document}

\vspace{-1in}

\maketitle
{\begin{abstract} Wilkie proved in 1977 that every countable model $\MM$ of Peano Arithmetic  has an elementary end extension $\NN$ such that the interstructure lattice $\Lt(\NN / \MM)$  is the pentagon lattice ${\mathbf N}_5$. This theorem implies that every countable nonstandard $\MM$ has an elementary cofinal extension $\NN$ such that $\Lt(\NN / \MM) \cong {\mathbf N}_5$. It is proved here that  whenever $\MM \prec \NN \models \pa$ and $\Lt(\NN / \MM) \cong{\mathbf N}_5$, then $\NN$ must be either an end or a cofinal extension of $\MM$. In contrast, there are $\MM^* \prec \NN^* \models \pas$ such that $\Lt(\NN^* / \MM^*) \cong{\mathbf N}_5$ and $\NN^*$ is neither an end nor a cofinal extension of $\MM^*$. \end{abstract}

\bigskip

Throughout, the   (possibly adorned) script letters $\MM, \NN, \KK$  denote models of Peano Arithmetic (\pa) having universes denoted by the (similarly adorned) roman letters $M,N,K$, respectively.  When we write $\MM \prec \NN$, we allow the possibility that $\MM = \NN$.  As usual,  we write $\MM \prec_{\sf end} \NN$ if $\NN$ is an {\em end} elementary extension of $\MM$ (that is, $a < b$ whenever $a \in M$ and $b \in N \backslash M$), and we  write $\MM \prec_{\sf cf} \NN$ if $\NN$ is a {\em cofinal} (necessarily  elementary) extension of $\MM$ (that is, for every $b \in N$ there is $a \in M$ such that $b < a$). If the elementary extension is neither end nor cofinal, then we say that it is {\bf mixed} and 
write $\MM \prec_{\sf mix} \NN$.

For a model ${\NN}$, its {\em substructure lattice} $\Lt(\NN)$   is the lattice of all those  $\KK \prec \NN$ ordered by $\prec$.  More generally, if $\MM \prec \NN$, then the {\em interstructure lattice} $\Lt(\NN / \MM)$ is  the sublattice of $\Lt(\NN)$ consisting of those ${\mathcal K}$ in $\Lt(\NN)$ 
such that $\MM \prec \KK$.   The question of which finite lattices can be  substructure (or, equivalently, interstructure, by Corollary~1.2) lattices is discussed in \cite[Chap.~4]{ksbook}. 
It is still unknown whether there are any finite lattices that are not substructure lattices; however, many lattices are known to be, among which  are   all the finite distributive lattices. In fact (\cite[Cor.~4.3.8]{ksbook}), for any for any $\MM$ and any finite distributive lattice $D$, there is $\NN \succ_{\sf end} \MM$ such that $\Lt(\NN / \MM) \cong D$.  

Recall that a lattice is distributive iff it embeds neither the pentagon lattice ${\mathbf N}_5$ nor the diamond lattice ${\mathbf M}_3$, both of which are depicted in Fig.~\ref{hex}. A lattice is modular iff it does not embed ${\mathbf N}_5$. 

\begin{figure}[h]
\begin{picture}(470,90)(-75,30)

%\put(-8,27){\line(1,0){240}}
%\put(232,27){\line(0,1){86}}
%\put(-8,113){\line(1,0){240}}
%\put(-8,27){\line(0,1){86}}
%\put(-11,24){\line(1,0){246}}
%\put(235,24){\line(0,1){92}}
%\put(-11,116){\line(1,0){246}}
%\put(-11,24){\line(0,1){92}}

% ----------------------------------

%\put(47,56){\circle{6}}
%\put(47,84){\circle{6}}
%\put(35,100){\circle{6}}
%\put(20,70){\circle{6}}
%\put(35,40){\circle{6}}

%\put(37,42){\line(3,4){8.6}}
%\put(33.4,42.5){\line(-1,2){12.3}}
%\put(33.4,97.5){\line(-1,-2){12.3}}
%\put(47,59){\line(0,1){22}}
%\put(37,98){\line(3,-4){8.6}}

%\put(27,105){\makebox(0,0){$1$}}
%\put(56,52){\makebox(0,0){$a$}}
%\put(56,87){\makebox(0,0){$b$}}
%\put(11,69){\makebox(0,0){$c$}}
%\put(27,37){\makebox(0,0){$0$}}

% ----------------------------------

%\put(157,40){\circle{6}}
%\put(157,80){\circle{6}}
%\put(137,60){\circle{6}}
%\put(157,60){\circle{6}}
%\put(177,60){\circle{6}}
%\put(157,63){\line(0,1){13}}
%\put(157,43){\line(0,1){13}}
%\put(159,41){\line(1,1){16}}
%\put(139,61){\line(1,1){16}}
%\put(155,41){\line(-1,1){16}}
%\put(159,78){\line(1,-1){16}}

% ----------------------------------

\put(47,56){\circle*{6}}
\put(47,84){\circle*{6}}
\put(35,100){\circle*{6}}
\put(20,70){\circle*{6}}
\put(35,40){\circle*{6}}

%\put(40,26){\makebox(0,0){$\nu_0$}}

\put(37,42){\line(3,4){8.6}}
\put(33.4,42.5){\line(-1,2){12.3}}
\put(33.4,97.5){\line(-1,-2){12.3}}
\put(47,59){\line(0,1){22}}
\put(37,98){\line(3,-4){8.6}}

\put(27,105){\makebox(0,0){$1$}}
\put(56,52){\makebox(0,0){$a$}}
\put(56,87){\makebox(0,0){$b$}}
\put(11,69){\makebox(0,0){$c$}}
\put(27,37){\makebox(0,0){$0$}}

% ----------------------------------

%\put(164,106){\makebox(0,0){$1$}}
%\put(122,78){\makebox(0,0){$a$}}
%\put(150,78){\makebox(0,0){$b$}}
%\put(191,78){\makebox(0,0){$c$}}
%\put(165,37){\makebox(0,0){$0$}}

\put(157,40){\circle*{6}}
\put(157,100){\circle*{6}}
\put(127,70){\circle*{6}}
\put(157,70){\circle*{6}}
\put(187,70){\circle*{6}}
\put(157,73){\line(0,1){23}}
\put(157,43){\line(0,1){23}}
\put(159,41){\line(1,1){26}}
\put(129,71){\line(1,1){26}}
\put(155,41){\line(-1,1){26}}
\put(159,98){\line(1,-1){26}}

\end{picture}
\caption{Lattices ${\mathbf N}_5$ and ${\mathbf M}_3$}
\label{hex}
\end{figure}

Paris \cite{par77} gave, historically, the first example of a substructure lattice that is not distributive.  
The following theorem of Wilkie \cite{wil77} gives the first example of a  substructure lattice that is not modular.

\bigskip

  {\sc Theorem 1}: {\em For every countable $\MM$ there is $\NN \succ_{\sf end} \MM$ such that $\Lt(\NN / \MM) \cong {\mathbf N}_5$.}
  
  \bigskip
  
  Incidentally, as proved in \cite[Th.~4.5]{ksbook}, for every  $\MM_0$ there is  $\MM \succ_{\sf end} \MM_0$  for which no $\NN \succ_{\sf end} \MM$ is such that $\Lt(\NN / \MM) \cong {\mathbf N}_5$. Theorem~1 has the following corollary. (See   Theorem~1.1 for the reason.)
  
   \bigskip
  
   {\sc Corollary 2}: {\em For every countable and nonstandard $\MM$ there is $\NN \succ_{\sf cf} \MM$ such that $\Lt(\NN / \MM) \cong {\mathbf N}_5$.}
   
    \bigskip
  
  It is still unresolved if, for every nonstandard $\MM$,  there is  $\NN \succ_{\sf cf} \MM$ such that $\Lt(\NN / \MM) \cong {\mathbf N}_5$. A positive answer  would immediately  yield a positive answer to the question
  \cite[Question~2, Chap.~12]{ksbook} if every uncountable model has a minimal cofinal extension.

 Theorem~1 and Corollary~2 together  suggest the question of whether the  pentagon lattice can be realized by an elementary mixed extension. It was impetuously stated in \cite[page~123]{ksbook} that ${\mathbf N}_5$ does have such a representation. There was no published proof at that time,    but there was an outline of a proof  that later was seen  to be flawed.   Mea culpa! In fact, there is no such extension.   That is the content of the next theorem,  which is our first main result.

 \bigskip
 
 {\sc Theorem 3}: {\em If $\MM \prec_{\sf mix} \NN$, then   $\Lt(\NN / \MM) \not\cong {\mathbf N}_5$.}
 
 \bigskip

 Despite Theorem~3, there is a sliver of truth to the claim in \cite{ksbook}, as will now be explained. 
 Let ${\mathcal L}$ be one of the usual finite languages in which  \PA\ is formulated; to be definitive, let ${\mathcal L} = \{+,\times,\leq,0,1\}$.
 Let ${\mathcal L}^*$ be the language obtained from  ${\mathcal L}$ by adjoining to it 
the  denumerably many new and distinct unary relation symbols $U_0,U_1,U_2, \ldots$. Thus, 
${\mathcal L}^*= {\mathcal L} \cup \{U_i : i < \omega\}$. 
Let $\pas$ be the ${\mathcal L}^*$-theory of  those structures $\MM^* = (\MM,U_0,U_1,U_2, \ldots)$, where $\MM \models 
\pa$ and $\MM^*$ satisfies the induction scheme  for all   ${\mathcal L}^*$-formulas. (From now on, $\MM^*,\NN^*, \ldots$ always denote models of $\pas$ that are expansions of $\MM, \NN, \ldots$.) We can think of $\pas$ as a subtheory of $\pa$ by identifying $\pa$ with $\pas \cup\{U_i = \varnothing : i < \omega\}$.  Many concepts, such as interstructure lattices,  that concern models of $\pa$ extend in an obvious and natural way to models of $\pas$. Also, 
many results  about models of \pa, together with their proofs,  extend in a straightforward manner to models of $\pas$. Almost all results in \cite{ksbook} do. Theorem~1 and Corollary~2 also do. In other words, Theorem~1$^*$ and Corollary~2$^*$ are valid, where we are adjoining  $^*$   to indicate that $\pa^*$ rather than just $\pa$ is being considered.
But Theorem~3 has the unusual feature that it does not. The next theorem, our  second new  result, indicates why.

\bigskip

\hyphenation{expan-sion}

{\sc Theorem 4}: {\em Every  countable, recursively saturated $\MM$ has an expansion  $\MM^*$ for which there is $\NN^* \succ_{\sf mix} \MM^*$ such that $\Lt(\NN^*/ \MM^*) \cong {\mathbf N}_5$.} 

\bigskip

 As far as notation and terminology go, we generally follow what is standard or what can be found in \cite{ksbook}.
 
 There are four numbered sections following this introduction. The first contains some preliminary material much of which is   a rehash. The second is almost purely combinatorial in nature and prepares  the way for  the proofs of Theorems~3 and~4, which are then presented in Sections~3 and~4, respectively.

 % ======================================================
 % ======================================================           
%  
%                               SECTION 1:  PRELIMINARIES
%
% ======================================================
% ======================================================

\section{Ranked Lattices and their Representations}  \label{sect1} This section,  comprising three subsections,  culminates with a description  of how to obtain elementary extensions realizing a given finite ranked lattice. The first subsection repeats some material from \cite[Chap.~4]{ksbook}. The second subsection extends the first and puts a new perspective on it. Finally, the third subsection   extends   the second from lattices to ranked lattices.

 \subsection{Representations of lattices} \label{sect11}For any set $A$, let $\Eq(A)$ be the lattice of equivalence relations on $A$,
ordered in such a way that if $\Theta_1,\Theta_2 \in \Eq(A)$,
then $\Theta_1 \leq \Theta_2$ iff $\Theta_1 \subseteq \Theta_2$ (that is,  $\Theta_1$ {\bf refines} $\Theta_2$).
We let  $\bbzero_A$ be the {\bf discrete} equivalence relation on  $A$ (that is, $\bbzero_A$ is the equality relation on $A$) and 
$\bbone_A$ be the {\bf trivial} equivalence relation (that is, $\bbone_A = A \times A$). Thus, 
for any $\Theta_1, \Theta_2 \in \Eq(A)$, we have that
$$
\bbzero_A \leq \Theta_1 \leq \bbone_A 
$$
and
$$
\Theta_1 \wedge \Theta_2 =  \Theta_1 \cap  \Theta_2. $$
If $\Theta \in \Eq(A)$ and $B \subseteq A$ then 
 $\Theta \cap B^2 \in \Eq(B)$. 
 If $f$ is a function with domain $A$, then $f$ {\bf induces} 
$\Theta \in \Eq(A)$
if whenever $a,b \in A$, then $\langle a , b \rangle \in \Theta$ iff 
$f(a) = f(b)$.

Let $L$ be a finite lattice.  
A {\bf representation} of $L$ is a one-to-one function 
$\alpha : L \into \Eq(A)$ such that 
$$
\alpha(0_L) = \bbone_A,
$$
$$
\alpha(1_L) = \bbzero_A,
$$
and 
$$
\alpha(r \vee s) = \alpha(r) \wedge  \alpha(s)
$$
 for each $r,s \in L$. (It is not required that
$\alpha(r \wedge s ) = \alpha(r) \vee \alpha(s)$.) We say that $\alpha$ is {\bf finite} if $A$ is a finite set.
   If $B \subseteq A$, then $\alpha|B : L \into \Eq(B)$ is  such that
$(\alpha|B)(r) = \alpha(r) \cap B^2$ for each $r \in L$.  The representation  $\beta : L \into \Eq(B)$ is {\bf isomorphic} to $\alpha$ (in symbols: $\alpha \cong \beta$) if there is a 
bijection $f : A \into B$ such that for any $x,y \in A$ and $r \in L$,
$\langle x , y \rangle \in \alpha(r)$ iff $\langle f(x) , f(y)\rangle \in  \beta(r)$. If such is the case, then we say that~$f$ {\bf demonstrates} that $\alpha \cong \beta$. 
If $\alpha : L \into \eq(A)$ is a representation and $\Theta \in \eq(B)$, then $\Theta$ is {\bf canonical} ({\bf for} $\alpha$) if  $B \subseteq A$ and  $\Theta = \alpha(r) \cap B^2$ for some $r \in L$. 

Suppose that $\alpha : L \into \eq(A)$ is a representation of the finite lattice $L$, and $\mathcal{B}$ is a set of representations of $L$. Then 
 $\alpha$ {\bf arrows} ${\mathcal B}$ (in symbols: $\alpha \into {\mathcal B}$) if whenever $\Theta \in \eq(A)$, then there is $B \subseteq A$ such that $\Theta \cap B^2$ is canonical for $\alpha$ and $\alpha|B \cong \beta$ for some $\beta \in {\mathcal B}$.   We usually write $\alpha \into \beta$ instead of $\alpha \into \{\beta\}$.
 
 We next define, by recursion on $n < \omega$, when the  representation $\alpha : L \into \eq(A)$ of the finite lattice $L$ has the 
 $n$-{\bf canonical partition property} (or, briefly, is $n$-{\bf CPP}).  First, $\alpha$ is $0$-CPP if for every $r \in L$, there do not exist exactly $2$ $\alpha(r)$-classes; next, $\alpha$ is $(n+1)$-CPP if there is a set ${\mathcal B}$ of $n$-CPP representations of $L$ such that $\alpha \into {\mathcal B}$. 
 
 Given $\MM$, we say that  a representation $\alpha$ of a finite lattice $L$ is an $\MM$-{\bf representation} 
 if it is $\MM$-definable. Also, if $A \in \Def(\MM)$, then we let $\eq^\MM(A)$ be the set of those $\Theta \in \eq(\MM)$ that are definable in $\MM$. All the definitions in this subsection up to this point make sense when interpreted in a model $\MM$ and are applied  just to $\MM$-representations. In particular, it makes sense to refer to an $\MM$-finite $\MM$-representation $\alpha$ as being $a$-CPP for $a \in M$. Thus, for every finite lattice $L$, there is a $\Sigma_1$ formula $cpp_L(x)$ such that for any $\MM$ and $a \in M$, $\MM \models cpp_L(a)$ iff there is an $\MM$-finite $\MM$-representation of $L$ that $\MM$ thinks is $a$-CPP. The following theorem can be found in \cite[Chap.~4]{ksbook} or \cite{sch86}.
 
 \bigskip
 
 {\sc Theorem 1.1}: {\em Let $L$ be a finite lattice and $\MM$ be a nonstandard, countable model. The following are equivalent$:$
 
 $(1)$ There are $\NN_0 \succ \MM_0 \equiv \MM$ such that  $\Lt(\NN_0 / \MM_0) \cong L$.
 
 $(2)$ For every $n < \omega$, $\MM \models cpp_L(n)$.
 
 $(3)$ There is $\NN \succ_{\sf cf} \MM$ such that $\Lt(\NN / \MM) \cong L$.} \qed
 
 \bigskip
 
 Notice that Corollary~2  follows from Theorem~1 and $(1) \Longrightarrow (3)$ of the previous theorem.
 
 \bigskip
 
 {\sc Corollary 1.2}: {\em If $L$ is a finite lattice, then the following are equivalent:
 
 \begin{itemize}
 
 \item[(1)] There is $\NN$ such that $\Lt(\NN) \cong L$.
 
 \item[(2)] There are $\MM \prec \NN$ such that $\Lt(\NN / \MM) \cong L$.
 
 \end{itemize}}
 
 \bigskip
 
 {\it Proof}. Obviously, $(1) \Longrightarrow (2)$ by letting $\MM$ be the prime elementary submodel of $\NN$. The converse  $(2) \Longrightarrow (1)$ follows from Theorem~1.1 as long as $\MM$ is not a model of True Arithmetic (\ta) and so its prime elementary submodel is nonstandard. If $\MM$ is a model of \ta, then by $(1) \implies (2)$ of Theorem~1.1, $\MM \models cpp_L(n)$ for all $n < \omega$. Since \ta\ is undecidable,  there is a prime, nonstandard  $\MM_0$ such that $\MM_0 \models  cpp_L(n)$ for all $n < \omega$. Then by $(1) \implies (2)$ of Theorem~1.1, there is $\NN_0 \succ_{\sf cf} \MM_0$ such that $\Lt(\NN_0) = \Lt(\NN_0 / \MM_0) \cong L$. \qed

%   ====> Section 1.2 : Correct sets   <========

\subsection{Correct Sets of Representations} \label{sect12} This subsection consists of  a definition followed by  a theorem generalizing 
 Theorem~1.1.
 
  \bigskip

{\sc Definition}  1.3: Let $\MM$ be a model 
and  $L$  be a finite  lattice. 
We say that ${\mathcal C}$ is an $\MM$-{\bf correct set of 
representations of} $L$ if each of the following hold.

\begin{itemize}

 \item[(1)] ${\mathcal C}$ is a nonempty set of $0$-CPP $\MM$-representations of $L$.

 \item[(2)] Whenever $\alpha : L \into \eq(A)$ is in ${\mathcal C}$ and 
$\Theta \in \Eq^\MM(A)$, 
then there is a $B \subseteq A$ such that $\alpha|B \in {\mathcal C}$
and $\Theta \cap B^2$ is canonical for $\alpha$.

\end{itemize}

 \bigskip
 
 Here is an example. Suppose that $\MM$ is nonstandard and that $\MM \models cpp_L(n)$ for every $n < \omega$. Let ${\mathcal C}$ be the set of those $\MM$-finite $\MM$-representations  $\alpha$ of $L$ such that, for some nonstandard $n \in M$, $\MM$ thinks that $\alpha$ is $n$-CPP. Then, ${\mathcal C}$  is $\MM$-correct. With this example, we see that the following theorem generalizes  a good portion of Theorem~1.1.
 It is a consequence of Theorem~1.1 when $\MM$ is countable and nonstandard. 
 
 \bigskip

{\sc Theorem}  1.4:  {\em   Suppose that $\MM$ is a model and $L$ is a finite  lattice.

\begin{itemize}

\item[(1)]
If there is $\NN \succ \MM$ such that  $\Lt(\NN / \MM) \cong  L$, 
then there is an $\MM$-correct set of
representations of $L$.

\item[$(2)$] If $\MM$ is  countable  and there is 
 an $\MM$-correct set of 
representations of $L$,
then there is $\NN \succ \MM$  such that $\Lt(\NN / \MM) \cong L$.

\end{itemize}} 

\bigskip

{\it Proof}. $(1)$ Suppose that $\NN \succ \MM$ and that $F : L \into \Lt(\NN / \MM)$ is an isomorphism.
Let $f : L \into N$ be such that for $r \in L$, $f(r)$ generates $F(r)$ over $\MM$. Let $a = f(1_L)$. 
 
 For each pair of elements $r,s \in L$, let $g_{r,s} : N \into N$ and $h_{r,s} : N^2 \into N$ be functions that are $\NN$-definable using parameters only from $M$ such that 
 \begin{itemize}
 
 \item $g_{r,s}(f(r \vee s)) = f(r)$,
 
 \item $h_{r,s}(f(r),f(s)) = f(r \vee s)$.
  
 \end{itemize}
 The functions $g_{r,s}$ exist since $f(r) \in F(r \vee s)$;  the functions $h_{r,s}$ exist  for  a similar reason.
 Let $g_r = g_{r,1}$, so that $g_r(a) = f(r)$. In particular, $g_1(a) = a$. The two equalities above become
 
 \begin{itemize}
 
 \item $g_{r,s}(g_{r \vee s}(a)) = g_r(a)$,
 
 \item $h_{r,s}(g_r(a),g_s(a)) = g_{r \vee s}(a)$.
 
 \end{itemize}
 
 For each $X \in \Def(\MM)$, let $\alpha_X : L \into \eq(X)$ be such that 
 whenever $r \in L$, then $\alpha_X(r)$ is the equivalence relation in $\eq(X)$ induced by 
 $g_r \harp X$.  Let $B$ be the set of all $x \in M$ such that 
  \begin{itemize}
 
 \item $g_{r,s}(g_{r \vee s}(x)) = g_r(x)$,
 
 \item $h_{r,s}(g_r(x),g_s(x)) = g_{r \vee s}(x)$,
 
 \item $g_1(x) = x$.
 
 \end{itemize} 
 Clearly, $B \in \Def(\MM)$ and $a \in B^\NN$. We claim that $\alpha_B$ is an $\MM$-representation of $L$. But even more is true. If $X \subseteq B$, $X \in \Def(\MM)$ and $a \in X^\NN$, then  $\alpha_X = \alpha_B|X$. 
  
 We now claim that each such $\alpha_X$ is an $\MM$-representation of $L$. 
 
 First, $\alpha_X$ is one-to-one. For, suppose that $r,s \in L$ and $\alpha_X(r) = \alpha_X(s)$. Then,
 $g_r \harp X$ and $g_1 \harp X$ induce the same equivalence relations on $X$. It follows that  there are $\MM$-definable functions $e_0,e_1 : M \into M$ such that for all $x \in X$, $e_0(g_r(x)) = g_s(x)$ and 
 $e_1(g_s(x)) = g_r(x)$.  But then $e_0^\NN(g_r(a)) = g_s(a)$ and 
 $e_1^\NN(g_s(a)) = g_r(a)$, implying that $F(r) = F(s)$ and, therefore, $r = s$.  
 
 Next, to prove that each $\alpha_X$ is a representation of $L$, it is enough to show that $\alpha_B$ is. 
 
 For all $x \in B$, $g_0(x) = f(0)$ and $g_1(x) = x$, so $\alpha_B(0)$ is trivial and $\alpha_B(1)$ is discrete. Finally, we show that if $r,s \in L$, then $\alpha_B(r \vee s) = \alpha_B(r)
\wedge \alpha_B(s)$. To do so, we let $x,y \in X$, and then show that 
$\langle x,y \rangle \in \alpha_B(r \vee s) \Longleftrightarrow \langle x,y \rangle \in \alpha_B(r) \cap \alpha_B(s)$.

\smallskip

\begin{eqnarray}
\langle x,y \rangle \in \alpha_B(r \vee s) & \Rightarrow & g_{r \vee s}(x) = g_{r \vee s}(y) \nonumber \\
 & \Rightarrow & g_{r,s}(g_{r \vee s}(x)) =  g_{r,s}(g_{r \vee s}(y))  \nonumber \\
& \Rightarrow &
g_r(x) = g_r(y)  \nonumber  \\
& \Rightarrow & \langle x,y \rangle \in \alpha_B(r). \nonumber
\end{eqnarray}
Similarly, $\langle x,y \rangle \in \alpha_B(r \vee s)  \Rightarrow \langle x,y \rangle \in \alpha_B(s)$.
Conversely, 

\begin{eqnarray}
\langle x,y \rangle \in \alpha_B(r) \cap \alpha_B(s) & \Rightarrow & g_r(x) = g_r(y) \ \& \ 
g_s(x) = g_s(y) \nonumber \\
& \Rightarrow & h_{r,s} (g_r(x),g_s(x)) = h_{r,s} (g_r(y),g_s(y)) \nonumber \\
& \Rightarrow & g_{r\vee s}(x) = g_{r\vee s}(y) \nonumber \\
& \Rightarrow & \langle x,y \rangle \in \alpha_B(r \vee s). \nonumber 
\end{eqnarray}

Having that each $\alpha_X$ is a representation of $L$, we easily see that it is $0$-CPP. 
For if $X$ is partitioned into  $Y,Z \in \Def(\MM)$, then either $a \in Y^\NN$ or $a \in Z^\NN$, but not both.

 Now let ${\mathcal C}$ be the set of all such $\alpha_X$; that is,
 $$ 
 {\mathcal C} = \{\alpha_X : X \subseteq B, X \in \Def(\MM),  a \in X^\NN\}.
 $$ 
 We have just seen that ${\mathcal C}$ is a nonempty set of $0$-CPP $\MM$-representations of $L$, so that $(1)$ of Definition~1.3 is verified. We prove (2) of Definition~1.3. Consider $\alpha_X \in {\mathcal C}$. Let $\Theta \in \eq(X)$ be  $\MM$-definable. Define $m : X \into X$ so that if $x \in X$, then $m = \min([x]_\Theta)$. 
 Let $r \in L$ be such that $m^\NN(a)$ generates $F(r)$ over $\MM$. There are functions 
 $e_0,e_1 : N \into N$ that are $\NN$-definable but using  parameters only from $M$  such that $e_0(m^\NN(a)) = r$ and $e_1(r) = m^\NN(a)$. Let $Y = \{x \in X :  e_0(m^\NN(a)) = r$ and $e_1(r) = m^\NN(a)\}$. 
 Then $m \harp Y$ induces $\alpha_Y(r)$. This completes the proof of (1).

\smallskip

$(2)$  Since ${\mathcal C} \neq \varnothing$, let $\alpha: L \into \eq(A)$ be in ${\mathcal C}$.  Let $\Theta_0,\Theta_1,\Theta_2, \ldots$ enumerate all $\MM$-definable equivalence relations on $M$. By recursion, obtain a sequence $X_0 \supseteq X_1 \supseteq X_2 \supseteq \cdots$ of sets in $\Def(\MM)$ 
   as follows. Let $X_0 = A$. Suppose that we have $X_n$ and that $\alpha|X_n \in {\mathcal C}$. Let $X_{n+1} \subseteq X_n$ be such that $\alpha|X_{n+1} \in {\mathcal C}$ 
and $\Theta_n \cap X_{n+1}^2$ is canonical for $\alpha$. The $X_n$'s generate a complete type over $\MM$ (using that each $\alpha|X_n$ is 0-CPP). Let $\NN$ be an elementary extension of $\MM$ generated by an element $a$ realizing this type.

For each $r \in L$, let $t_r : M \into M$ be $\MM$-definable such that 
whenever $x \in X_0$, then $t_r(x) = \min([x]_{\alpha(r)})$. Define the function $F$ on $L$ so that if $r \in L$, then  $F(r)$ is the elementary substructure of $\NN$ generated by $t^\NN_r(a)$ over $\MM$. One easily checks  that $F : L \into \Lt(\NN / \MM)$ is an isomorphism. \qed

\bigskip

\subsection{Ranked Lattices} \label{sect13}To  refine the notions of end/cofinal/mixed extensions, we appeal to rankings of lattices \cite[Def.\@ 4.2.6]{ksbook}. Suppose that $L$ is a  finite lattice. A function $\rho : L \into L$ is a {\bf ranking} of $L$ if for each $r,s \in L$:

\begin{itemize}
\item[(1)] $\rho(r) \geq r,$
\item[(2)] $\rho(\rho(r)) = \rho(r)$,
\item[(3)] $\rho(r) \leq \rho(s)$ or $\rho(s) \leq \rho(r)$,
\item[(4)] $\rho(r \vee s) = \rho(r) \vee \rho(s)$.

\end{itemize}
A ranking $\rho$ of $L$ uniquely determines, and is uniquely determined by, its {\bf rankset} 
$\{\rho(r) : r \in L\}$. If $L$ is finite and $R \subseteq L$, then $R$  is a rankset iff $R$ is linearly ordered 
and $1_L \in R$.
If $\rho$ is a ranking of $L$, then $(L,\rho)$ is a {\bf ranked lattice}. 

If $\MM \prec \NN$ and $\Lt(\NN / \MM)$ is finite, then let $\rho : \Lt(\NN / \MM) \into \Lt(\NN / \MM)$ be such that if $\KK  \in  \Lt(\NN / \MM)$, then $\rho(\KK)$ 
 is uniquely defined by
$$
\KK \prec_{\sf cf} \rho(\KK) \prec_{\sf end} \NN.
$$
One easily verifies that $\rho$ is a ranking of $ \Lt(\NN / \MM)$. We let $\Ltr(\NN / \MM)$ be the ranked lattice $(\Lt(\NN / \MM),\rho)$.

 Suppose that $\rho$ is a ranking of the finite lattice $L$. Then  $\rho$ 
is an {\bf end} ranking if $\rho(0_L) = 0_L$, a  {\bf cofinal} ranking if $\rho(0_L) = 1_L$ and a {\bf mixed} ranking if $0_L < \rho(0_L) < 1_L$.  Obviously, $L$ has a unique cofinal ranking. If $\rho$ is an end,  cofinal or mixed ranking, then $(L,\rho)$ is, respectively, an {\bf end},  {\bf cofinal} or {\bf mixed} ranked lattice. These definitions are  appropriate:\@
if $\Ltr(\NN / \MM)$  is an end,  cofinal or mixed  ranked lattice, then   $\NN$ is, respectively,  an end,  cofinal or mixed extension of $\MM$.

\begin{figure}[h]
\begin{picture}(470,90)(-75,30)

%\put(-8,27){\line(1,0){240}}
%\put(232,27){\line(0,1){86}}
%\put(-8,113){\line(1,0){240}}
%\put(-8,27){\line(0,1){86}}
%\put(-11,24){\line(1,0){246}}
%\put(235,24){\line(0,1){92}}
%\put(-11,116){\line(1,0){246}}
%\put(-11,24){\line(0,1){92}}

% ----------------------------------

\put(47,56){\circle{6}}
\put(47,84){\circle{6}}
\put(35,100){\circle*{7}}
\put(20,70){\circle{6}}
\put(35,40){\circle{6}}

\put(40,26){\makebox(0,0){$\nu_0$}}

\put(37,42){\line(3,4){8.6}}
\put(33.4,42.5){\line(-1,2){12.3}}
\put(33.4,97.5){\line(-1,-2){12.3}}
\put(47,59){\line(0,1){22}}
\put(37,98){\line(3,-4){8.6}}

% ----------------------------------

\put(107,56){\circle{6}}
\put(107,84){\circle*{7}}
\put(95,100){\circle*{7}}
\put(80,70){\circle{6}}
\put(95,40){\circle*{7}}

\put(100,26){\makebox(0,0){$\nu_1$}}

\put(97,42){\line(3,4){8.6}}
\put(93.4,42.5){\line(-1,2){12.3}}
\put(93.4,97.5){\line(-1,-2){12.3}}
\put(107,59){\line(0,1){22}}
\put(97,98){\line(3,-4){8.6}}

% ---------------------------------

\put(167,56){\circle*{7}}
\put(167,84){\circle*{7}}
\put(155,100){\circle*{7}}
\put(140,70){\circle{6}}
\put(155,40){\circle*{7}}

\put(160,26){\makebox(0,0){$\nu_2$}}

\put(157,42){\line(3,4){8.6}}
\put(153.4,42.5){\line(-1,2){12.3}}
\put(153.4,97.5){\line(-1,-2){12.3}}
\put(167,59){\line(0,1){22}}
\put(157,98){\line(3,-4){8.6}}

% -------------------------------------

\put(227,56){\circle*{7}}
\put(227,84){\circle*{7}}
\put(215,100){\circle*{7}}
\put(200,70){\circle{6}}
\put(215,40){\circle{6}}

\put(220,26){\makebox(0,0){$\nu_3$}}

\put(217,42){\line(3,4){8.6}}
\put(213.4,42.5){\line(-1,2){12.3}}
\put(213.4,97.5){\line(-1,-2){12.3}}
\put(227,59){\line(0,1){22}}
\put(217,98){\line(3,-4){8.6}}

\end{picture}
\caption{Four Ranked Pentagon Lattices}
\label{ranked}
\end{figure}

Of the  $10$ rankings of ${\bf N}_5$, four  are depicted in Fig.~\ref{ranked} by letting 
\ $\bullet$ \ denote  those points in  the rankset and
\ $\circ$ \ those that are not. Of all the ranked pentagons, the four  in Fig.~\ref{ranked} are 
the most important for us because of the following. 

Henceforth, we use the labeling of ${\mathbf N}_5$ as given in Fig.~\ref{hex}.

\bigskip

{\sc Proposition 1.5}: {\em If $\MM \prec \NN$ and $\Ltr(\NN / \MM) \cong ({\mathbf N}_5, \rho)$, then 
$\rho = \nu_i$ for some $i \leq 3$.}

  \bigskip
  
  {\it Proof}.  We first show that $\rho(c) = 1$.  If $\rho(c) \neq 1$, then $\rho(c) = c$. We apply the Gaifman Condition \cite[Prop.\@ 4.2.12]{ksbook} by letting $x = a$, $y = b$ and $z =c$,  to get the contradiction  that $a = b$. 
  
  If $\rho(0) = 1$, then $\rho = \nu_0$. So, assume that $\rho(0) < 1$. 
   Since $\rho(c) = 1$ and $c \wedge b = 0$, it follows from the Blass Condition \cite[Prop.\@ 4.2.7]{ksbook} that $\rho(b) = b$. Finally, $\rho(0) \neq b$ by \cite[Thm.\@ 4.6.1]{ksbook}. Thus, $\rho(0) \in \{0,a\}$, so it must be that $\rho \in \{\nu_1,\nu_2, \nu_3\}$. \qed
   
   \bigskip
   
   We make some comments about this proposition. First, Proposition~1.5$^*$ is also valid.  Theorem~1 can now be restated as: For all countable $\MM$ there are $i \in \{1,2\}$ and $\NN \succ \MM$  such that  
$\Ltr(\NN / \MM) \cong({\mathbf N}_5, \nu_i)$. In fact, Wilkie's proof of Theorem~1 yields that $i = 1$. A similar proof shows that for every countable $\MM$ there is $\NN \succ_{\sf end} \MM$ such that $\Ltr(\NN / \MM) \cong ({\mathbf N}_5, \nu_2)$. Since $\nu_3$ is the only mixed ranking of the four in  Fig.~\ref{ranked}, then in Theorem~4 we get 
$\NN^*$ such that $\Ltr(\NN^* / \MM^*) \cong ({\mathbf N}_5, \nu_3)$.

\smallskip

The next order of business is to generalize Definition~1.3 and Theorem~1.4  from lattices to ranked lattices. 

First, we need some terminology. Suppose that $\MM$ is a model, $A \in \Def(\MM)$ and 
$\Theta \in \Eq^\MM(A)$. 
We say that a set ${\mathcal E}$ of $\Theta$-classes is  $\MM$-{\bf bounded} if there is a bounded $I \in \Def(\MM)$  such that $I \cap X \neq \varnothing$ for each $X \in {\mathcal E}$.

If $(L,\rho)$ is a finite ranked lattice, then a representation $\alpha$ of $L$   is a {\bf representation of} $(L,\rho)$ if whenever $r \leq s \in L$, then 
$s \leq \rho(r)$ iff every $\alpha(r)$-class is the union of a finite set of $\alpha(s)$-classes. This definition  should help motivate  the next definition.

  \bigskip
 
 {\sc Definition 1.6}: Let  $\MM$ be a model and $(L,\rho)$  a finite ranked lattice. 
 
(1)  A representation $\alpha : L \into \eq(A)$  is an  $\MM$-{\bf representation} of $(L,\rho)$ if    $\alpha$ is an   
 $\MM$-representation of $L$  and  whenever $r \leq s \in L$, then 
$s \leq \rho(r)$ iff every $\alpha(r)$-class is the union of an $\MM$-bounded set of $\alpha(s)$-classes.

  (2)
  We say that ${\mathcal C}$ is an $\MM$-{\bf correct set of 
representations of} $(L,\rho)$ if  ${\mathcal C}$ is an $\MM$-correct set of representations of $L$ and 
 each $\alpha \in {\mathcal C}$ is an $\MM$-representation of $(L,\rho)$. 

\bigskip

We next generalize Theorem~1.4 from lattices to ranked lattices.

 \bigskip

{\sc Theorem}  1.7:  {\em   Suppose that $\MM$ is a model and $(L,\rho)$ is a finite  ranked lattice.

\begin{itemize}

\item[(1)]
If there is $\NN \succ \MM$ such that  $\Ltr(\NN / \MM) \cong  (L,\rho)$, 
then there is an $\MM$-correct set of
representations of $(L,\rho)$.

\item[$(2)$] If $\MM$ is  countable  and there is 
 an $\MM$-correct set of 
representations of $(L,\rho)$,
then there is $\NN \succ \MM$  such that $\Ltr(\NN / \MM) \cong (L,\rho)$.

\end{itemize}} 

\bigskip

 {\it Proof}. $(1)$ Obtain ${\mathcal C}$ as in the proof of Theorem~1.4(1), so that ${\mathcal C}$ is an 
$\MM$-correct set of representations of $L$.  If $\alpha : L \into\eq(A)$ is in  ${\mathcal C}$, $r \leq s$ but not $s \leq \rho(r)$, then there is some $\alpha(r)$-class that is not the union of an $\MM$-bounded set of $\alpha(s)$-classes. (For, otherwise, there would be an $\MM$-definable function $b : M \into M$ such that $b(f(r)) \geq f(s)$.) However, it could be that  $r \leq s \leq \rho(r)$ and some $\alpha(r)$-class is not the union of an $\MM$-bounded set of $\alpha(s)$-classes. Let ${\mathcal C}_0$ be the set of those $\alpha \in {\mathcal C}$ that are $\MM$-representations of $(L, \rho)$. We will show that this ${\mathcal C}_0$ is  an $\MM$-correct set of representations of $(L,\rho)$. To see this, it suffices to show that for each $\alpha : L \into \eq(A)$  in ${\mathcal C}$, there is $B \subseteq A$ such that $\alpha|B \in {\mathcal C}_0$.

Suppose that we have $\alpha : L \into \eq(A)$  in ${\mathcal C}$ and that $r \leq s \leq \rho(r)$.
Partition $A$ into two sets $A_0,A_1$, so that $A_0$ is the union of those $\alpha(r)$-classes that are the union of an $\MM$-bounded set of $\alpha(s)$-classes. Since ${\mathcal C}$ is $\MM$-correct, then either $\alpha|A_0 \in {\mathcal C}$ or $\alpha|A_1 \in {\mathcal C}$. By what was previously said, the latter option is impossible, so we have that $\alpha|A_0 \in {\mathcal C}$. Repeating this for all such $r,s \in L$, finally yields $B \subseteq A$ as required. This completes the proof of  (1).

\smallskip

(2)  Let ${\mathcal C}$ be an $\MM$-correct set of representations of $(L,\rho)$.
Then ${\mathcal C}$ is an $\MM$-correct set of representations of $L$, so we can obtain 
$\NN \succ \MM$ as in the proof of Theorem~1.4(2). Then $\Lt(\NN / \MM) \cong L$. 

We use the notation from the proof of Theorem~1.4(2). Thus, 
$F : L \into \Lt(\NN / \MM)$ is an isomorphism and $F(r)$ is generated by $t_r(a)$ over $\MM$. 
We prove that  $F$ is also an isomorphism of the ranked lattices 
$(L,\rho)$ and $\Ltr(\NN / \MM)$.   It suffices to prove: whenever $r < s \in L$, then $s \leq \rho(r)$ iff  
$F(r) \prec_{\sf cf} F(s)$. So, let $r < s \in L$.

\smallskip

$(\Longrightarrow)$: Suppose that $s \leq \rho(r)$. Consider $\alpha \in {\mathcal C}$. Every $\alpha(r)$-class is the union of an $\MM$-bounded set of $\alpha(s)$-classes. Let $g : M \into M$ be 
an $\MM$-definable function such if $x \in X$, then $g(x) = \max\{t_s(y) : \langle x,y \rangle \in \alpha(r)\}$. Clearly, $g(x)$ is well defined for $x \in X$, so there is such an $\MM$-definable $g$.
Thus, $g(t_r(x)) \geq t_s(x)$ for all $x \in X$, so that $g^\NN(t_r^\NN(a)) \geq t_s^\NN(a)$.
Therefore, $F(r) \prec_{\sf cf} F(s)$.                                                  

\smallskip

$(\Longleftarrow)$: Suppose that $F(r) \prec_{\sf cf} F(s)$. Then there is an $\MM$-definable 
$g : M \into M$ such that $g^\NN(t_r^\NN(a)) \geq t_s^\NN(a)$. 
There is $X_i$ such that $g(t^r(x)) \geq t_s(x)$ for all $x \in X_i$. Let $\alpha_i = \alpha|X_i \in {\mathcal C}$. 
Thus, each $\alpha_i(r)$-class is the union of an $\MM$-bounded set of $\alpha_i(s)$-classes. 
Then 
$s \leq \rho(s)$. \qed

\bigskip

Wilkie's proof of Theorem~1 made implicit use of Theorem~1.7(2).

% ==================================================================
%
%                             SECTION 2 : Representations of N_5  
%
% ==================================================================
\bigskip

\section{Representations of ${\mathbf N}_5$}  \label{sect2}For almost all of this section, we ignore  \pa\ and concentrate just on representations of ${\mathbf N}_5$. Only in the first and last paragraphs is \pa\ considered.

 \smallskip
 
 Caveat lector: In the next definition, and throughout this paper, $\omega^n$ is not an ordinal but is the set of $n$-tuples of natural numbers.  If $s \in \omega^n$ and $i < n$, then $s_i$ is the $i$-th element of $s$. Also, remember that if $n < \omega$, then $n = \{0,1, \ldots, n-1\}$.
 If $s \in \omega^n$ and $i < m \leq n$, then $s \harp m \in \omega^m$ and $(s \harp m)_i = s_i$.

 \bigskip
 
 {\sc Definition 2.1}: 
 For  $n < \omega$,  let $A_n = (n+2) \times \omega^{n+1}$ and then define
$\alpha_n : {\mathbf N}_5 \into \eq(A_n)$ so that  
 $\alpha_n(0)$ is trivial, $\alpha_n(1)$ is discrete, and whenever $i,j \leq n+1$ and $s,t \in \omega^{n+1}$, then 
 \begin{itemize} 
 \item $\big\langle \langle i,s \rangle, \langle j,t \rangle \big\rangle \in \alpha(a)$ iff $i = j$;
  
 \item $\big\langle \langle i,s \rangle, \langle j,t \rangle \big\rangle \in \alpha(b)$ iff $i = j$ and $s \harp i  = t \harp j$;
 
 \item $\big\langle \langle i,s \rangle, \langle j,t \rangle \big\rangle \in \alpha(c)$ iff $s = t$. 
 
 \end{itemize}
 
 \bigskip
 
 It is clear that each $\alpha_n$ is a representation of ${\mathbf N}_5$ and that, if   $n \geq 1$, then $\alpha_n$ is  $0$-CPP. Observe that if $m < n < \omega$ and $I \subseteq n+2$ is such that $|I| = m+2$, then there is $D \subseteq \omega^{n+1}$  such that $\alpha_m \cong \alpha_n|(I \times D)$. Our primary goal in this section is to prove the following theorem, which will be given a more precise formulation in Theorem~2.8.
 
 \bigskip
 
 {\sc Theorem 2.2}: {\em If $ m < \omega$, then there is $n < \omega$ such that $\alpha_n \into \alpha_m$.}
 
 \bigskip
 
 To prove this theorem, we will take a detour  and visit some other lattices and their representations. These  lattice are introduced in Definition~2.3 and their representations in Definition~2.4.    
    \bigskip
  
  {\sc Definition 2.3}: Suppose that $1 \leq m \leq n < \omega$.
   Let $G_{m,n}$ be the set consisting of all pairs $\langle \theta, f \rangle$, where $\theta \in \eq(n+1)$ and $f : n+1 \into m+1$ are such that if $i,j \leq n$ and $\langle i,j \rangle \in \theta$, then $f(i) = f(j)$. Let $\trieq$ be the partial ordering of $G_{m,n}$ such that if $\langle \theta, f \rangle, \langle \psi, g \rangle \in G_{m,n}$, then 
    $$\langle \theta, f \rangle \trieq \langle \psi, g \rangle {\mbox{ iff }} \theta \supseteq \psi {\mbox{ and }} f(i) \leq g(i) {\mbox{ for all }} i \leq n.
  $$ 
  \bigskip
  
  Clearly, $\trieq$ really is a partial ordering. It should be observed that $G_{m,n}$ with $\trieq$, as in Definition~2.3,  is a lattice in which
  $$
  0_{G_{m,n}} = \langle \bbone_{n+1},0 \rangle,
  $$
    $$
  1_{G_{m,n}} = \langle \bbzero_{n+1},m \rangle,
  $$
  $$
  \langle \theta,f \rangle \vee \langle \psi, g \rangle = \langle \theta \cap \psi,\sup(f,g) \rangle,
  $$
  where $\sup(f,g) = h$ iff $h(i) = \max(f(i),g(i))$ for all $i \leq n$.
  In the above equalities, we are identifying $k \leq m$ with the function  that is constantly $k$ on $n+1$. We will continue to do so.
  
  Our real concern is with the lattices $G_n = G_{n,n}$. The more general $G_{m,n}$ are introduced in order to be able to do an inductive proof.
  One of the reasons for introducing the lattices $G_n$ is that there is an embedding $e_n : {\mathbf N}_5 \into G_n$ defined by:
   \begin{eqnarray*}
  e_n(0) & = & 0_{G_n}, \\
  e_n(a) & =  & \langle \bbzero_{n+1}, 0 \rangle, \\
  e_n(b) & = & \langle \bbzero_{n+1}, \id_{n+1}\rangle, \\
  e_n(c) & = & \langle \bbone _{n+1},n \rangle, \\
  e_n(1) & = & 1_{G_n}.
  \end{eqnarray*}
   As usual, $\id_X$ is the identity function on $X$.

  \setlength{\unitlength}{.88pt}

\begin{figure}[h]
\begin{picture}(470,160)(-170,-15)

%\put(47,60){\makebox(0,0){$L$}}
%\put(73,60){\makebox(0,0){$K$}}

\put(32,133){\makebox(0,0){$1$}}
\put(57,18){\makebox(0,0){$a$}}
\put(57,104.5){\makebox(0,0){$b$}}
\put(-17,68){\makebox(0,0){$c$}}
\put(33,-9){\makebox(0,0){$0$}}
%\put(29,66){\makebox(0,0){$r$}}

%\put(63,104.5){\makebox(0,0){$d$}}

\put(20,-7.5){\circle*{8}}
\put(20,130){\circle*{8}}
%\put(49,61.3){\oval (35,60.3)}
%\put(72,60.5){\oval (35,70)}

\put(49,29){\circle*{8}}
\put(49,91.5){\circle*{8}}
\put(-4,60.75){\circle*{8}}
%\put(73,95.5){\circle*{8}}

\put(22,127){\line(3,-4){25}}
\put(18,127){\line(-1,-3){21}}
\put(23,-6){\line(3,4){24}}
\put(19,-5){\line(-1,3){21}}

\put(49,31){\line(0,1){62}}

\put(20,61.25){\circle*{8}}
\put(20,61.25){\line(0,1){70}}

\put(20,65.25){\line(3,-4){25}}
\end{picture}
\caption{Embedding ${\mathbf N}_5$ into $G_1$}

\label{G2}
\end{figure}

 \setlength{\unitlength}{.88pt}

It is routine to verify that each $e_n$ is an embedding. Fig.~3 depicts the lattice $G_1$ with ${\mathbf N}_5$ embedded in it. If $r \in {\mathbf N}_5$, then $e_1(r)$ is labeled with $r$. The unlabeled point is  $\langle \bbzero_2, 1-\id_2 \rangle$, where  $1-\id_2$ is the function $f : 2 \into 2$ such that 
$f(0) = 1$ and $f(1) = 0$.

 \bigskip
 
 Next, we define representations of the $G_{m,n}$.

  \bigskip
  
  {\sc Definition 2.4}: Suppose that $1 \leq m \leq n < \omega$.  Let $\gamma_{m,n} : G_{m.n} \into \eq\big((n+1) \times \omega^m\big)$ be such that if $\langle \theta, f \rangle \in G_{m,n}$ and $\langle i, s \rangle, \langle j, t \rangle \in (n+1) \times \omega^m$, then $\big\langle \langle i, s \rangle, \langle j, t \rangle \big\rangle \in \gamma_{m,n}( \langle \theta, f \rangle)$ iff 
  $\langle i,j \rangle \in \theta$ and   $s \harp f(i) = t \harp f(j)$. 
  
  \bigskip
  
    Observe that  $\gamma_{m,n}$ is indeed a representation of $G_{m,n}$.  However, no $\gamma_{m,n}$ is $0$-CPP since if $\theta$ has exactly 2 equivalence classes, then $\gamma_{m,n}(\langle \theta, 0 \rangle)$ has exactly 2 equivalence classes.  In fact, if $E$ is a $\theta$-class, then $E \times \omega^m$ is a $\gamma_{m,n}(\langle \theta, 0 \rangle)$-class. Thus, the number of $\gamma_{m,n}(\langle \theta, 0 \rangle)$-classes is equal to the number of $\theta$-classes. On the other hand, if $f : n+1 \into m+1$ is not constantly $0$, then there are infinitely many $\gamma_{m,n}(\langle \theta,f \rangle)$-classes.
  
   Let $\gamma_n = \gamma_{n,n}$. Note that for $ n < \omega$, both of the representations $\gamma_{n+1}$ and $\alpha_n$ are into $\eq\big((n+2) \times \omega^{n+1}\big)$. In fact, even more is true.
   
   \bigskip
   
   {\sc Lemma 2.5}: {\em If $n < \omega$, then  $\alpha_n = \gamma_{n+1} \circ e_{n+1}$.} 
   
   \bigskip
   
   {\it Proof}. The routine proof is left to the reader. \qed
   
   \bigskip
   
   We next come to the main result about the representations $\gamma_{m,n}$.

  \bigskip
  
  {\sc Lemma 2.6}: {\em If $1 \leq m \leq n < \omega$, then $\gamma_{m,n} \into \gamma_{m,n}$.}
  
  \bigskip

  {\it Proof}. To prove this lemma, we are given $1 \leq m \leq n < \omega$ and $\Theta \in \eq\big((n+1) \times \omega^m\big)$. The goal is to get $X \subseteq (n+1) \times \omega^m$ and $r \in G_{m,n}$ such that $\gamma_{m,n}|X \cong \gamma_{m,n}$ and  $\Theta \cap X^2 = \gamma_{m,n}(r) \cap X^2$. Observe that if $X \subseteq (n+1) \times \omega^m$ and $\gamma_{m,n}|X \cong \gamma_{m,n}$, then there is $D \subseteq \omega^m$ such that 
  $X = (n+1) \times D$.  This can be seen by observing that  each $\gamma_{m,n}(\langle \bbone_{n+1},m \rangle)$-class has the form $(n+1) \times \{s\}$ for some $s \in \omega^m$.  
  
  We first prove the special case of Lemma~2.6  in which the number of $\Theta$-classes is finite. 
  
  \bigskip
  
  {\sc Lemma 2.6.1} {\em If $1 \leq m \leq n < \omega$ and $\Theta \in \eq\big((n+1) \times \omega^m\big)$ has only finitely many equivalence classes, then there are $X \subseteq (n+1) \times \omega^m$ and $r \in G_{m,n}$ such that $\gamma_{m,n} | X \cong \gamma_{m,n}$ and $\Theta \cap X^2 = \gamma_{m,n}(r) \cap X^2$.}
  
  \bigskip
  
  {\it Proof}. The proof is by induction on $m \geq 1$; that is, for each $m$, we will prove the lemma  for all  $n \geq m$ and appropriate $\Theta$.
  
  Notice that the sought for $r \in G_{m,n}$  will necessarily be such that $\gamma_{m,n}(r)$ has only finitely many equivalence classes. This implies that $r = \langle \theta,0 \rangle$ for some $\theta \in \eq(n+1)$.  Recall that the $X$ we are after will be such that $X = (n+1) \times D$ for some $D \subseteq \omega^m$.

  \smallskip  
  
  {\sf  Basis step} $m = 1$: Consider some $\Theta \in \eq\big((n+1) \times \omega\big)$. Let $\Psi$ be the equivalence relation on $\omega$ such that 
  $\langle k,\ell \rangle \in \Psi$ iff for all $i \leq n$, $\big\langle \langle i,k \rangle, \langle i,\ell \rangle \big\rangle \in \Theta$. Since there are only finitely many $\Psi$-classes, there is an infinite $\Psi$-class $D$. Then, $X = (n+1) \times D$ is as required.
  Notice that $r = \langle \theta,0 \rangle$ for some $\theta \in \eq(n+1)$.
  
  \smallskip
  
  {\sf Inductive step} $m > 1$: We are assuming that the lemma has been proved for all smaller values of $m$. For each $s \in \omega^{m-1}$, let $\Psi_s$ be the equivalence relation on $\omega$ such that $\langle k, \ell \rangle \in \Psi_s$ iff there are $t,t' \in \omega^m$ such that $t,t' \supseteq s$, $t_{m-1} = k$, $t'_{m-1} = \ell$ and $\big\langle \langle i,t \rangle, \langle i,t' \rangle \big\rangle \in \Theta$ for all $i \leq n$. Since there are only finitely many $\Psi_s$-classes, there is an infinite $\Psi_s$-class $Y_s$ such that $\Psi_s \cap Y_s^2$ is trivial. Now define $\Theta' \in \eq\big((n+1) \times \omega^{m-1}\big)$ so that 
  $\big\langle \langle i,s \rangle, \langle j,s' \rangle \big\rangle \in \Theta'$ iff for some (equivalently: all) 
  $t,t' \in \omega^m$ such that $t \supseteq s$, $t' \supseteq s'$, $t_{m-1} \in Y_s$ and $t'_{m-1} \in Y_{s'}$, 
  then $\big\langle \langle i,t \rangle, \langle j,t' \rangle \big\rangle \in \Theta$. There are only finitely many 
  $\Theta'$-classes, so by the inductive hypothesis, there are $\theta' \in \eq(n+1)$ and  $D' \subseteq \omega^{m-1}$ such that, letting $Y = (n+1) \times D'$, then $\gamma_{m-1,n} | Y \cong \gamma_{m-1,n}$ and $\Theta' \cap Y^2 = \gamma_{m-1,n}(\langle \theta' ,0 \rangle) \cap Y^2$. 
  Let 
  $$
  D = \{t \in \omega^m: {\mbox{ for some }} s \in D', t \supseteq s {\mbox{ and }} t_{m-1} \in Y_s\}. 
  $$
  It can be verified that there is $\theta \in \eq(n+1)$ such that $\theta' = \theta \cap (n \times n)$ and, letting $X = (n+1) \times D$, then $\gamma_{m,n} | X \cong \gamma_{m,n}$ and $\Theta \cap X^2 = \gamma_{m,n}(\langle \theta ,0 \rangle) \cap X^2$. \qed
  
   \bigskip

  {\it Proof of Lemma 2.6}. The proof is by induction on $m \geq 1$ with  $n \geq 1$ being fixed. The basis step is for $m = 1$ and the inductive step for $m > 1$. Both steps start out the same way. So for now, consider $m,n \geq 1$ and let $\Theta \in \eq\big((n+1) \times \omega^m\big)$.
  
   Consider an arbitrary $s \in \omega^{m-1}$. With the idea of invoking Infinite Ramsey's Theorem for pairs,  we define $F_s : [\omega]^2 \into \eq\big(\{0,1\} \times (n+1)\big)$ so that whenever $\{k,\ell \} \in [\omega]^2$, $k < \ell$, \ $e,e' \in \{0,1\}$ and $i,j \leq n$, then $\big\langle \langle e,i \rangle, \langle e',j \rangle \big\rangle \in F_s(\{k,\ell\})$ iff there are $t,t' \in \omega^m$ such that $t_{m-1} = k$, $t'_{m-1} = \ell$, $t \harp (m-1) = t' \harp (m-1) = s$ and one of the following:
   
   \begin{itemize}
   \item $e = e' = 0$ and $\big\langle \langle i,t \rangle, \langle j,t \rangle \big\rangle \in \Theta$;
   
    \item $e = 0$, $e' = 1$ and $\big\langle \langle i,t \rangle, \langle j,t' \rangle \big\rangle \in \Theta$;
    
      \item $e = 1$, $e' = 0$ and $\big\langle \langle i,t' \rangle, \langle j,t \rangle \big\rangle \in \Theta$;
      
        \item $e = e' = 1$ and $\big\langle \langle i,t' \rangle, \langle j,t' \rangle \big\rangle \in \Theta$.
        \end{itemize}
    Now we apply  Ramsey to get an infinite $H_s \subseteq \omega$ such that $F_s \harp [H_s]^2$ is constant. Let
   $$
   Y_s = \{t \in \omega^m :  t \supseteq s \mbox{ and } t_{m-1} \in H_s\}
   $$
   and
   $$
   X_0 = \bigcup \{(n+1) \times Y_s : s \in \omega^{m-1}\}.
   $$
   It is readily seen that $\gamma_{m,n} | X_0 \cong \gamma_{m,n}$  and that 
   each of the following is true for each $s \in \omega^{m-1}$:
   
   \begin{itemize}
   
   \item[(1)] If $i \leq n$, then $\Theta \cap (\{i\} \times Y_s)^2$ is either trivial or discrete.
   
   \item[(2)] If $i, j \leq n$, $\Theta \cap ( \{i\} \times Y_s)^2$ is trivial, $\Theta \cap  (\{j\} \times Y_s)^2$ is discrete and $t,t' \in Y_s$, then $\big\langle \langle i,t \rangle, \langle j,t' \rangle \big\rangle \not\in \Theta$.
   
   \item[(3)] If $i < j \leq n$ and both $\Theta \cap  (\{i\} \times Y_s)^2$ and $\Theta \cap  (\{j\} \times Y_s)^2$ are discrete,
   then one of the following: 
   
   \begin{itemize}
   \item[(3a)] if $t,t' \in Y_s$, then $\big\langle \langle i,t \rangle, \langle j,t' \rangle \big\rangle \not\in \Theta$;

   \item[(3b)] if $t,t' \in Y_s$, then $\big\langle \langle i,t \rangle, \langle j,t' \rangle \big\rangle \in \Theta$ iff $t = t'$.
   
   \end{itemize}
   \end{itemize}
   A consequence of (1) -- (3) is:
   
   \begin{itemize}
   \item[(4)] If $i,j \leq n$ and $t,t' \in Y_s$, then 
   $$
   \big\langle \langle i,t \rangle, \langle j,t \rangle \big\rangle \in \Theta \Longleftrightarrow \big\langle \langle i,t' \rangle, \langle j,t' \rangle \big\rangle \in \Theta.
   $$
   \end{itemize}
   Let 
   $T_s = \{i \leq n : \Theta \cap ( \{i\} \times Y_s)^2$ is trivial$\}$. Because of (1), 
   $(n+1) \backslash T_s = \{i \leq n : \Theta \cap ( \{i\} \times Y_s)^2$ is discrete$\}$.
   With (4) in mind, we can let $\theta_s \in \eq(n+1)$ be such that 
   $$
   \theta_s = \{\langle i,j \rangle \in (n+1) \times (n+1) : \big\langle \langle i,t \rangle, \langle j,t \rangle \big\rangle \in \Theta\}
   $$ 
   for each $t \in Y_s$.
   Clearly,  $T_s$ is the union of $\theta_s$-classes.
   
   \smallskip
   
   {\sf Basis step $m = 1$}: Since $m=1$,  there is only one possible $s$, namely $s = \varnothing$. 
   Thus, $X_0 = (n+1) \times Y_{\varnothing}$. We have already noted that $\gamma_{1,n}|X_0 \cong \gamma_{1,n}$. To complete this step, we need to show that there is $r_0 \in G_{1,n}$ such that $\Theta \cap X_0^2 = \gamma_{1,n}(r_0) \cap X_0^2$.

    Let $f : n+1 \into 2$ be such that $f(i) = 0$ iff $i \in T_{\varnothing}$.      Then we can take $r_0 = \langle \theta_\varnothing,f \rangle$. One easily verifies that  $\Theta \cap X_0^2 = \gamma_{1,n}(r_0) \cap X_0^2$.

   \smallskip
   
   {\sf Inductive step } $m > 1$: Thus, we are assuming  $\gamma_{m-1,n} \into \gamma_{m-1,n}$.  
   We already have infinite $Y_s$ for each $s \in \omega^{m-1}$ and  that (1) -- (4) hold.  Also, we have $X_0 \subseteq (n+1) \times \omega^m$ such that $\gamma_{m,n} | X_0 \cong \gamma_{m,n}$.  Without loss of generality, we  assume, for each $s \in \omega^{m-1}$, that  $H_s = \omega$ and then $Y_s = \{t \in \omega^m : t \supseteq s\}$. Thus, 
   $X_0 = (n+1) \times \omega^m$. 
      
   Let $\Theta_1 \in \eq\big((n+1) \times \omega^{m-1}\big)$ be such that  if $t \in Y_s$, $t' \in Y_{s'}$ and $i,j \leq n$, then $\big\langle \langle i,t \rangle, \langle j,t' \rangle \big\rangle \in \Theta_1$ iff $T_s = T_{s'}$ and $\theta_s = \theta_{s'}$. Since there are only finitely many $\Theta_1$-classes, we can apply  Lemma~2.6.1 to get $D'_1 \subseteq \omega^{m-1}$ and $X'_1 = (n+1) \times D'_1$ such that $\gamma_{m-1,n} | X'_1 \cong \gamma_{m-1,n}$ and $X'_1$ is a subset of some $\Theta_1$-class. Thus, there are $T$ and $\theta$ such that $T_s = T$ and $\theta_s = \theta$ whenever $s \in D'_1$. Let $D_1 = \{t \in \omega^m : t \harp (m-1) \in D'_1\}$ and let $X_1 = (n+1) \times D_1$. We then have that $\gamma_{m,n} | X_1 \cong \gamma_{m,n}$. Without loss of generality, we will assume that $D_1 = \omega^m$ so that $X_1 = (n+1) \times \omega^m$.  Notice that (1) -- (4) remain true and, in addition, the following hold:
   
   \begin{itemize}
   
   \item[(5)] If $i \leq n$, then $\Theta \cap (\{i\} \times Y_s)^2$ is trivial iff $i \in T$.
   
   \item[(6)] If $i,j \leq n$ and $t \in \omega^m$, then $\langle i,j \rangle \in \theta$ iff $\big\langle \langle i,t \rangle, \langle j,t \rangle \big\rangle \in \Theta$.
   
   \end{itemize}
   
    Let $\Theta_2 \in \eq\big((n+1) \times \omega^{m-1}\big)$ be such that  if $i,j \leq n$ and $s,s' \in \omega^{m-1}$, then $\big\langle \langle i,s \rangle, \langle j,s' \rangle \big\rangle \in \Theta_2$ iff one of the following:
   
   \begin{itemize}
   \item $i,j \not\in T$, $\langle i,j \rangle \in \theta$ and $s = s'$;
  
   \item $i,j \in T$ and 
   for some (or, equivalently, all) $t \in Y_s$ and $t'  \in Y_{s'}$, 
   $\big\langle \langle i,t \rangle, \langle j,t' \rangle \big\rangle \in \Theta$.
   \end{itemize}
   One easily verifies that, indeed,  $\Theta_2 \in \eq\big((n+1) \times \omega^{m-1}\big)$. By the inductive hypothesis, 
   there are $D'_2 \subseteq \omega^{m-1}$, $X'_2 =  (n+1) \times D_2'$ and $r_2' \in G_{m-1,n}$ such that 
   $\gamma_{m-1,n} | X'_2 \cong \gamma_{m-1,n}$ and $\Theta_2 \cap (X'_2)^2 = \gamma_{m-1,n}(r'_2) \cap (X'_2)^2$. Let $r'_2 = \langle \psi, f' \rangle$. It must be that $\psi = \theta$. 
   
   Now let $D_2 = \{t \in \omega^m : t \harp (m-1) \in D'_2\}$ and $X_2 = (n+1) \times D_2$.
   Then $\gamma_{m,n} | X_2 \cong \gamma_{m,n}$. Let $r_2 = \langle \psi,f \rangle \in G_{m,n}$, where
   $f(i) = f'(i)$ if $i \in T$ and $f(i) = m$ if $i \not\in T$. It may not be the case that 
   $\Theta \cap X^2_2 = \gamma_{m,n}(r_2) \cap X_2^2$, but we do have 
   $$
   \Theta \cap (T \times D_2)^2 = \gamma_{m,n}(r_2) \cap (T \times D_2)^2.
   $$
   Without loss of generality, let $D_2 = \omega^m$, so that $X_2 = (n+1) \times \omega^m$. Thus, we have, in addition to (1) -- (6), that:
   
   \begin{itemize}
   
   \item[(7)]  $\Theta \cap (T \times \omega^m)^2 = \gamma_{m,n}(r_2) \cap (T \times \omega^m)^2.$
   
   \end{itemize}
   
   To complete this inductive step, we proceed with what might be called a ``thinning'' of $\omega^m$. The object is to get $D_3 \subseteq \omega^m$ such that if $X_3 = (n+1) \times D_3$, then 
   $\gamma_{m,n} | X_3 \cong \gamma_{m,n}$ and $\gamma_{m,n}(r_2) \cap X^2_3 = \Theta \cap X^2_3$.
   
   Let $\langle s^k : k < \omega \rangle$ be a one-to-one enumeration of $\omega^m$. By recursion on $k$, choose  
   $t^k \in \omega^m$ so that:
   \begin{itemize}
   \item[(T1)] $t_k \not\in \{t^0,t_1, \ldots,t^{k-1}\}$,
   
   \item[(T2)] $t^k \harp (m-1) = s^k \harp (m-1)$,
   
   \item[(T3)] $\Theta$ and $\gamma_{m,n}(r_2)$ agree on $(n+1) \times \{t_0,t_1, \ldots, t^k\}$.
   
   \end{itemize}
   
   Clearly, $t^0 = s^0$.  If $k > 0$, then there are only finitely many $t \in \omega$ such that (T2) holds but (T3) fails, so it is always possible to get $t^k$. 
   
   Let $r = r_2$, $D = \{t^k : k < \omega\}$ and $X = (n+1) \times D$. Then, $X$ and $r$ are as required. \qed

\bigskip

{\sc Corollary 2.7}: {\em If $1 \leq n < \omega$, then $\gamma_n \into \gamma_n$.}

\bigskip

    {\it Proof}. Let $n = m$ in Lemma~2.6. \qed
    
    \bigskip
    
    Let $R : \omega \into \omega$ be the Ramsey function such that if $m < \omega$, then  $R(m)$  is the least $k < \omega$ such that whenever  $\chi : [k]^2 \into 35$, then there is $I \subseteq k$ such that $|I| = m$ and $\chi$ is constant on $[I]^2$. (It seems to be right that 35 is large enough for the following proof to work. But if it isn't, replace it with something that is.)
    
    \bigskip
    
    {\sc Theorem} 2.8: {\em If $m < \omega$ and $n \geq 4R(m+2)^2$, then $\alpha_n \into \alpha_m$.}
    
    \bigskip
    
    {\it Proof}. Let $\Theta \in \eq\big((n+2) \times \omega^{n+1}\big)$. By Lemma~2.5 and Corollary~2.7, we can assume that $\Theta = \gamma_n(\langle \theta,f \rangle)$, where $\langle \theta,f \rangle \in G_{n+1}$. 
    
    Since $\theta \in \eq(n+2)$, there is $J \subseteq n+2$ such that $|J| = 2R(m+2)$ and $\theta \cap J^2$ is either trivial or discrete. We consider each of these possibilities.
    
    \smallskip
    
    {\it trivial}: The function $f \harp J$ is constantly $r \leq n+1$. Since $R(m+2) \geq m+2$,  there is $I \subseteq J$ such that 
    $|I| = m+2$ and either (1) $i > r$ for each $i \in I$ or (2) $i < r$ for each $i  \in r$. There is $D \subseteq \omega^{n+1}$ such that  $\alpha_n |(I \times D) \cong \alpha_m$. Let $X = I \times D$. If (1), then 
    $\Theta \cap X^2 = \alpha_n(0)$; and if (2), then $\Theta \cap X^2 = \alpha_n(c)$.
    
    \smallskip
    
    {\it discrete}: Let $\chi : [J]^2 \into 35$ be such that if $i,i',j,j' \in J$, $i < j$, $i' < j'$ and $\chi(\{i,j\}) = \chi(\{i',j'\})$, then $\{\langle i,i' \rangle, \langle j,j' \rangle, \langle f(i),f(i') \rangle, \langle f(j),f(j') \rangle\}$ 
     is an order-preserving function.  Let $I \subseteq J$ be such that $|I| = m+2$ and $\chi$ is constant on $[I]^2$. There are 3 possibilities: (1) $f(i) < j$ for all $i,j \in I$; (2) $f(i) = i$ for all $i \in I$; (3) $f(i) > j$ for all $i,j \in I$. There is $D \subseteq \omega^{n+1}$ such that 
     $\alpha_n |(I \times D) \cong \alpha_m$. Let $X = I \times D$. If (1), then 
    $\Theta \cap X^2 = \alpha_n(a)$;  if (2), then $\Theta \cap X^2 = \alpha_n(b)$; and if (3), then $\Theta \cap X^2 = \alpha_n(1)$.
    
    \smallskip
    This completes the proof. \qed

 \bigskip

  A careful inspection of the previous proofs shows that  they can be carried out in $\aca$. (Keep in mind that $(\MM, \Def(\MM)) \models \aca$ for every $\MM$, where $\Def(\MM)$ is the set of definable subsets of $M$.) For example, we see from the proof of Lemma~2.6 that if $1 \leq m \leq n < \omega$, then $\aca \vdash \gamma_{m,n} \into \gamma_{m,n}$. The function $R$ defined just before Theorem~2.8 can defined in any model of $\pa$, and then also $\aca$. Thus, we get the following theorem.

  \bigskip
  
  \vspace{40pt} \qed
   
   \vspace{-40pt}

  {\sc Theorem 2.9}: {\em If $m < \omega$ and $n = 4R(m+2)^2$, then 
  $$
   \aca \vdash \alpha_n \into  \alpha_m. 
   $$}

  \bigskip

  % =================================================== 
  %                                                                                                             
  %   ======>   SECT. 3 : Proving THEOREM 3  <============          
  %                                                                                                         
  % ==================================================
  
  \section{Proving Theorem~3 } \label{sect3} This section is devoted to proving Theorem~3.    
 With the idea of obtaining a contradiction, assume  that $\MM \prec_{\sf mix} \NN$ and that 
$\Lt(\NN / \MM) \cong {\bf N}_5$.  Proposition~1.4 implies that $\Ltr(\NN / \MM)$ $ \cong ({\bf N}_5, \nu_3)$. Following Theorem~1.7(1),  we let ${\mathcal C}$ be an $\MM$-correct set of representations of $({\mathbf N}_5, \nu_3)$.  In the course of this proof, we will see that ${\mathcal C}$ must have certain properties. We will also see that there are other properties that ${\mathcal C}$ possibly could have, and we will then assume that ${\mathcal C}$ does have these properties. 

Since ${\mathcal C} \neq \varnothing$, fix some $\alpha \in {\mathcal C}$.   Thus, $\alpha : {\mathbf N}_5 \into \eq(A)$. We can assume: 

\begin{itemize}

\item[(C1)] {\em For every $\beta \in {\mathcal C}$, there is $B \subseteq A$ such that $\beta = \alpha | B$.}

\end{itemize}

Since $a \vee c = 1$, then $\alpha(a) \cap \alpha(c) = \bbzero_A$; therefore, whenever  $X$ is an $\alpha(a)$-class 
  and $Z$ is an $\alpha(c)$-class,  then $|X \cap Z| \leq 1$. Since $0 < a = \nu_3(0)$, then, according to Definition~1.6, 
   the set of $\alpha(a)$-classes is $\MM$-bounded; let  $n +1\in M$ be the number of $\alpha(a)$-classes according to $\MM$. Then, we can assume:    
   \begin{itemize}
   
   \item[(C2)] $\alpha : {\mathbf N}_5 \into \eq(A)$, where   $A = [0,n] \times M$ and $n$ is nonstandard, is such that if $\langle i,j \rangle, \langle i',j' \rangle \in A$, then
   $$
   \big\langle \langle i,j \rangle, \langle i',j' \rangle \big\rangle \in \alpha(a) \mbox{ iff } i = i',
   $$
   $$
   \big\langle \langle i,j \rangle, \langle i',j' \rangle \big\rangle \in \alpha(c) \mbox{ iff } j = j'.
   $$
   
   \end{itemize}
   At first, it may look as if we can only assume that $A \subseteq [0,n] \times M$. But it is always possible to enlarge the set $A$ so as to get $[0,n] \times M$.
   
   \smallskip

   For just this proof, let us say that the $\MM$-representation $\beta$ of  ${\mathbf N}_5$  is {\it rectangular} if  $|X \cap Z| = 1$ for each $\beta(a)$-class $X$ and $\beta(c)$-class $Z$.  We see from (C2) that $\alpha$ is rectangular. We can even assume:
   
   \begin{itemize}
   \item[(C3)] 
      {\em Every $\beta \in {\mathcal C}$ is a  rectangular representation.} 
      \end{itemize}      
     
     To see why, 
let ${\mathcal C}_0$ be the set of those rectangular $\MM$-representations $\beta$, where $B \subseteq A$ and $\beta = \alpha |B$, 
   for which  there is $A_0 \subseteq B$ such that    $\alpha|A_0 \in {\mathcal C}$. To prove that this ${\mathcal C}_0$ is an $\MM$-correct set of representations of $({\mathbf N}_5,\nu_3)$, it suffices to prove 
   that if $A_1 \subseteq A$ and $\alpha | A_1\in {\mathcal C}$, then there is  $B \subseteq A_1$ such that $\alpha|B \in {\mathcal C}_0$. To prove this, 
   consider some $\alpha_1 = \alpha|A_1 \in {\mathcal C}$.  Define $\Theta \subseteq \eq(A_1)$ so 
   that if $y,z \in A_1$, then $\langle y,z \rangle \in \Theta$ iff the following holds for  each  
   $\alpha_1(a)$-class $X$: there is $u \in X$ such that $\langle u,y \rangle \in  \alpha_1(c)$ iff 
   there is $v \in X$ such that $\langle v,z \rangle \in \alpha_1(c)$. Clearly, $\alpha_1(c) \subseteq  \Theta \in  \eq(A_1)$.  Since ${\mathcal C}$ is $\MM$-correct,
   there are $A_0 \subseteq A_1$ and $r \in \{0,c\}$ such that $\alpha|A_0 \in {\mathcal C}$ and 
   $\alpha_1(r) \cap A_0^2 = \Theta \cap A_0^2$.  The number of $\Theta$-classes is at most $2^{n+1}$,    so it must be that $r = 0$. Let $B$ be the union of those $\alpha_1(c)$-classes  
   that have a nonempty intersection with $A_0$.  Then $A_0 \subseteq B \subseteq A_1$ and   $\beta = \alpha|B \in {\mathcal C}_0$.  This proves that ${\mathcal C}_0$ is an $\MM$-correct set of rectangular representations of $({\mathbf N}_5, \nu_3)$, so we can assume (C3).
   
   \smallskip
   
   Moreover, we can also assume:
   
      \begin{itemize}
   \item[(C4)] 
      {\em If $I \subseteq I' \subseteq [0,n]$, $J \subseteq J' \subseteq M$ and $\alpha|(I \times J)  \in {\mathcal C}$, then $\alpha|(I' \times J') \in {\mathcal C}$.} 
      \end{itemize}

Working in $\MM$, let $\langle B_k : k \in M \rangle$ be a  one-to-one enumeration of all $\alpha(b)$-classes. Thus, there is an ${\mathcal L}(M)$-formula $\psi(u,v)$ such that 
 $$
 \MM \models \forall u,v[\psi(u,v) \leftrightarrow v \in B_u].
 $$
 If $\beta : {\mathbf N}_5 \into \eq(X)$ is in ${\mathcal C}$ and $p \in M$, then there is $X' \subseteq X$ 
 such that $\beta | X' \in {\mathcal C}$ and if $B_k \cap X' \neq \varnothing$, then $k > p$. To see why, just consider $\Theta \in \eq(X)$ such that $\bigcup\{B_k \cap X : k > p\}$ is a $\Theta$-class and then apply Definitions~1.3 and~1.6.  

For each $j \in M$, there is a (unique) permutation $\pi_j$ of $n$ defined by the following condition: 
if $i,i' < n$, then $\pi_j(i) \leq \pi_j(i')$ iff there are (necessarily unique) $k,k' \in M$ such that 
$\langle i,j \rangle \in B_{k}, \langle i',j \rangle \in B_{k'}$ and $k \leq k'$. Using these permutations, 
we define $\Psi \in \eq(A)$ so that $\big\langle \langle i,j \rangle, \langle i',j' \rangle \big\rangle \in \Psi$ iff $\pi_{j} = \pi_{j'}$. Clearly, the set of $\Psi$-classes is $\MM$-bounded  as $\MM$ thinks that there are  no more than  $(n+1)!$ $\Psi$-classes. 
Thus, there are $I \times J \subseteq [0,n] \times M = A$ and $\pi$ such that $\alpha|(I \times J) \in {\mathcal C}$ and $\pi_j = \pi$ whenever $\langle i,j \rangle \in I \times J$. 
Without loss of generality, we  assume that $J = M$ and that  $\pi$ is the identity 
permutation. Thus, we have:

\begin{itemize}
\item[(C5)] If $i,i' \leq n$, $j \in M$, $\langle i,j \rangle \in B_k$ and $\langle i',j \rangle \in B_{k'}$, then
$$
i \leq i' \mbox{ iff } k \leq k'.
$$
\end{itemize} 

We now have that ${\mathcal C}$ is an $\MM$-correct set  of representations  of $({\mathbf N}_5,\nu_3)$ satisfying (C1) -- (C5).

With (C5) in mind, we make a couple of definitions concerning a $\beta \in {\mathcal C}$, where $\beta : {\mathbf N}_5 \into \eq(I \times J)$.  Suppose that $X$ and $Y$ are $\beta(b)$-classes. We say that $X$ is {\it below} $Y$ if there are $i,i' \in I$ and $j \in J$ such that $\langle i,j \rangle \in X$, $\langle i',j \rangle \in Y$ and $i < i'$. If $X$ is below $Y$, then $Y$ is {\it above} $X$. Thus, (C4) says: If $B_k$ is below $B_{k'}$ (as $\alpha(b)$-classes), then $k < k'$. 
The following is a consequence of (C5):

\begin{itemize}

\item[(C6)] For each $\alpha(b)$-class, the set of $\alpha(b)$-classes below it is $\MM$-bounded. 

\end{itemize}

\bigskip

{\sc Lemma} 3.1: {\em For each $\beta \in {\mathcal C}$, there is a $\beta(b)$-class $X$ having an $\MM$-unbounded set of $\beta(b)$-classes above it.}

\bigskip

{\it Proof}. Suppose that $\beta : {\mathbf N}_5 \into \eq(I \times J)$ is in ${\mathcal C}$. Let $\Theta \in \eq(I \times J)$ be such that whenever $\langle i,j \rangle, \langle i',j' \rangle \in I \times J$, then 
$\big\langle \langle i,j \rangle, \langle i',j' \rangle \big\rangle \in \Theta$ iff for every $\beta(b)$-class $X \subseteq \{i''\} \times J$, $\langle i'',j \rangle \in X$ iff $\langle i'',j' \rangle \in X$. Clearly,
$\beta(c) \subseteq \Theta \in \eq^\MM(I \times J)$. Thus, there are $I' \times J' \subseteq I \times J$ and $r \in \{0,c\}$ such that $\beta|(I' \times J') \in {\mathcal C}$ and $\Theta \cap (I' \times J')^2 = \beta(b) \cap (I' \times J')^2$.

If $r = 0$, then every $(\beta(b) \cap (I' \times J')^2)$-class is a  $(\beta(a) \cap (I' \times J')^2)$-class, which is impossible.
Hence, $r = c$. 

Let $X'$ be an $\MM$-unbounded $(\beta|(I' \times J'))(b)$-class. Then there  is an $\MM$-unbounded set of $(\beta|(I' \times J'))(b)$-classes above it. If $X$ is the $\beta(b)$-class such that $X \supseteq X'$, then the set of $\beta(b)$-classes above $X$ is $\MM$-unbounded. \qed

\bigskip

We define the function $q : [0,n] \times M \into M$ so that if $\langle i,j \rangle \in [0,n] \times M$, then 
$q(i,j)$ is such that $\langle i,j \rangle \in B_{q(i,j)}$. 

\bigskip

{\sc Corollary 3.2}: {\em If $\beta : {\mathbf N}_5 \into \eq(I \times J)$ is in ${\mathcal C}$ and $p \in M$, then there is $X \subseteq I \times J$ such that $\beta|X \in {\mathcal C}$ and 
if $\langle i,j \rangle \in X$, then   $q(i,j)  > p$.}

\bigskip

*********************

\smallskip

Consider some $\beta : {\mathbf N}_5 \into \eq(B)$ in ${\mathcal C}$, where $B = I \times J \subseteq [0,n] \times M$. 
For each $m < \omega$, we say that $\beta$ is $m$-{\it thick} if whenever $f : M \into M$ is definable, then there are $i_0,i_1, \ldots,i_m \in I$ and $j \in J$ such that $f(q(i_k,j)) < q(i_{k+1},j)$ for each $k < m$.

\bigskip

{\sc Lemma} 3.3: {\em If $ m < \omega$, then every  $\beta \in {\mathcal C}$ is $m$-thick.}

\bigskip

{\it Proof}. The proof is by induction on $m$. Notice that every $\beta  \in {\mathcal C}$ is vacuously $0$-thick.

\smallskip

{\sf Basis step} m = 1: Let $\beta : {\mathbf N}_5 \into \eq(B)$ be in ${\mathcal C}$. By Lemma~3.1, there is a $\beta(b)$-class $B \cap B_{k_0}$ having an $\MM$-unbounded set of $\beta(b)$-classes above it. Thus, for any $\MM$-definable $f : M \into M$, there is a $\beta(b)$-class  $B \cap B_{k_1}$ above $B \cap B_{k_0}$ such that $f(k_0) < k_1$. 

\smallskip

{\sf Inductive step} $2 \leq m < \omega$: We are  assuming that every $\beta \in {\mathcal C}$ is $(m-1)$-thick, where $\beta : {\mathbf N}_5 \into \eq(I \times J)$. Consider some definable $f : M \into M$. 
We define $g : J \into (I \cup \{\max(I) + 1\})$. For each $j \in J$, let $g(j) \in I$ be the least such that there are $i_0,i_1, \ldots,i_{m-1} = g(j)$ such that $f(q(i_k,j)) < q(i_{k+1},j)$ for all $k < m-1$. (If there is no such $g(j)$, then let $g(j) = \max(I) +1$.) By the inductive hypothesis, there is $j \in J$ such that $g(j) \in I$. Notice that $g$ is $\MM$-definable.

We define $\Theta \in \eq(I \times J)$ having exactly  2 $\Theta$-classes which are 
$D_0 = \{ \langle i,j \rangle \in I \times J : i < g(j)\}$ and $D_1 = (I\times J) \backslash D_0$. 
Since $\Theta$ is $\MM$-definable, there is $I' \times J' \subseteq I \times J$ such that 
$\beta|(I' \times J') \in {\mathcal C}$ and either $I' \times J' \subset D_0$ or $I' \times J' \subset D_1$. 
The first alternative contradicts the inductive hypothesis, so we have that $I' \times J' \subset D_1$. 

By Lemma~3.1, there are $i,i',j$ such that $\langle i,j \rangle, \langle i' , j \rangle \in I' \times J'$ and $f(q(i,j)) < q(i',j)$. There are $i_0,i_1, \ldots i_{m-2} \in I$ such that 
$f(q(i_k,j)) < q(i_{k+1},j)$ for $k < m-2$ and $f(q(i_{m-2},j) < q(i,j)$. Let $i_{m-1} = i$ and $i_m = i'$.
Then, $i_0,i_1, \ldots,i_m$ and $j$ demonstrate that $\beta$ is $m$-thick. \qed

  \bigskip 
  
  Lemma~3.3 implies a strengthening of itself via Corollary~3.2.
  
  \bigskip
  
  {\sc Corollary 3.4}: {\em  Suppose that $m < \omega$, $p \in M$, $\beta : {\mathbf N}_5 \into \eq(X)$ in ${\mathcal C}$, where $X = I \times J$, and $f : M \into M$ is definable. Then  there are $i_0,i_1, \ldots,i_m \in I$ and $j \in J$  such that $q(i_0,j) > p$ and $f(q(i_k,j)) < q(i_{k+1},j)$ for each $k < m$. } \qed
  
  \bigskip
  
  We need some more notation and terminology. Recall that a  cut $K$ (of $\MM$) is a  subset of $M$ such that $0 \in K \neq M$  and that $x + 1 \in K$ whenever $x \leq y \in K$. If $m < \omega$, then the cut $K$ is  $\Sigma_m$-{\it closed} iff whenever $\varphi(x)$ is a $\Sigma_m$ ${\mathcal L}(K)$-formula and $\MM \models \exists x \varphi(x)$, then there is $a \in K$ such that $\MM \models \varphi(a)$.  If $K$ is a $\Sigma_0$-closed cut and $\varphi \in {\mathcal L}(K)$ is a formula, then   $\ulcorner \varphi \urcorner$, the G\"odel number of $\varphi$, is in $K$.
  
  Also, recall that $\alpha : {\mathbf N}_5 \into \eq(A)$, $A = [0,n] \times M$ and $\langle B_k : k \in M \rangle$ is a definable, one-to-one enumeration of the $\alpha(b)$-classes (as defined a few lines after (C3)).  Working in $\MM$, let $\Gamma$ be the set of all ${\mathcal L}(M)$-sentences that are in prenex form with prenex ${\mathsf Q}_0 x_0 {\mathsf Q}_1 x_1   \cdots  {\mathsf Q}_i x_i   \cdots{\mathsf Q}_{n} x_n$, where  each ${\mathsf Q}_i$ is either $\forall$ or $\exists$ and whose matrix has length at most $n$.  (Any standard sentence can be thought of as being in $\Gamma$.) For each $\sigma \in \Gamma$ and $j \in M$, let $\sigma^{(j)}$ be the sentence derived from $\sigma$ by replacing each $Q_ix_i$ with $Q_ix_i \leq q(i,j)$. Thus, $\MM$ thinks that each $\sigma^{(j)}$ is a $\Delta_0$ sentence, so the set of those $\sigma^{(j)}$, where $\sigma \in \Gamma$ and $j \in M$, that $\MM$ thinks are true is definable in $\MM$. 
  
   In $\MM$, define the function $F$  on $A$ so that if $\langle i,j \rangle \in A$, then $F(i,j)$ is the set of all 
   ${\mathcal L}([0,q(i,j)])$-sentences $\sigma \in \Gamma$ such that
  $\MM$ thinks that $\sigma^{(j)}$ is true. Observe that if $\langle i,j \rangle, \langle i',j' \rangle \in A$ and $F(i,j) = F(i',j')$, then $q(i,j) = q(i',j')$.
 Let $\Theta \in \eq(A)$ be the equivalence relation induced by $F$; that is, if $\langle i,j \rangle, \langle i', j' \rangle \in A$, then $\big\langle \langle i,j \rangle, \langle i', j' \rangle \big\rangle \in \Theta$ iff $F(i,j) = F(i',j')$.   Thus, we have that $\Theta \in \eq^\MM(A)$ and  $\Theta \subseteq \alpha(b)$.   By Definition~1.6, there is $X = I \times J \subseteq A$ such that 
 $\alpha | X \in {\mathcal C}$ and $\Theta \cap X^2 = \alpha(r) \cap X^2$ for some $r \in {\mathbf N}_5$.
 Obviously, $r \in \{b,1\}$. But since, $F[B]$ is $\MM$-bounded for each $\alpha(b)$-class $B$, it must be that $r = b$.   Let $\beta = \alpha | X$. 
  
  Choose a sufficiently large $m < \omega$ and a sufficiently large $d \in M$. (Here, by ``sufficiently large'', we mean that $d > n^n$ and there is  $k$ such that $k+2 < m$ and each of $\alpha, X,q$ is definable by a $\Sigma_k$ $ {\mathcal L}([0,d])$-formula.)
  
   Let $g : M \into M$ be a one-to-one  function that is defined by a $\Sigma_{m+1}$ ${\mathcal L}$-formula that is fast-growing in the sense that for every $f : M \into M$ that is defined by $\Sigma_m$ ${\mathcal L}$-formula and every nonstandard $x \in M$, $\MM \models g(x) > f(x)$.    
   By Corollary~3.4, we can let $i_0,i_1, \ldots, i_m \in I$ and $j \in J$ be such that $q(i_0,j) > g(d)$ and 
   $g(q(i_\ell,j)) < q(i_{\ell+1})$ for each $\ell < m$. Let $K$ be the smallest $\Sigma_0$-closed cut such that $d \in K$. Thus, we have that for each standard $\Sigma_m$ ${\mathcal L}(K)$ sentence $\sigma$, $\MM \models \sigma \leftrightarrow \sigma^{(j)}$. But also, for every $j' \in J$ such that $q(i_0,j') = q(i_0,j)$, we have that $\MM \models \sigma^{(j)} \leftrightarrow \sigma^{(j')}$.  Then there is a $\Sigma_{m-1}$ ${\mathcal L}(K)$-formula defining the true $\Sigma_m$ ${\mathcal L}(K)$-sentences.   But this contradicts  the following version of Tarski's Theorem on the undefinability of truth, which is an immediate consequence of G\"odel's Diagonalization Lemma.

 \bigskip
 
 {\sc Theorem 3.5}: {\em Suppose that $1 \leq m < \omega$, $K \subseteq M$ is a $\Sigma_0$-closed cut and $\varphi(x)$ is an ${\mathcal L}(K)$-formula such that for each $\Sigma_m$ ${\mathcal L}(K)$-sentence 
 $\sigma$, $\MM \models \varphi(\ulcorner \sigma \urcorner) \leftrightarrow \sigma$. Then $\varphi(x)$ is not a $\Pi_m$ formula.} \qed
 
 \bigskip
 
 This contradiction completes the proof of Theorem~3.

\bigskip

  % =================================================== 
  %                                                                                                             
  %   ======>   SECT. 4 : Proving THEOREM 4  <============          
  %                                                                                                         
  % ==================================================

\section{Proving Theorem~4} \label{sect4} This section is devoted to a proof of Theorem~4. The definitions of ${\mathcal L}^*$ and the ${\mathcal L}^*$-theory $\pas$ are 
given in the introduction.  For each $\MM^*$, observe  that $(\MM, \Def(\MM^*)) \models \aca$. For any countable, recursively saturated 
$\MM$ we will obtain $\MM^*$ and $\NN^*$ as in that theorem. First, 
 we isolate a certain class of models of $\pas$ that have such extensions. 

\bigskip

{\sc Definition 4.1}: A model $\MM^*$ is {\bf recursively supersaturated} if 
$(\MM, \Def(\MM^*))$ (qua a two-sorted, first-order model of second-order arithmetic)
is recursively saturated.

\bigskip

{\sc Proposition 4.2}: {\em Every countable, recursively saturated $\MM$ can be expanded to a recursively supersaturated $\MM^*$.}

\bigskip

{\it Proof.} Let $T = \Th(\MM) + \aca$. Then $T \in \ssy(\MM)$, so it has a countable $\ssy(\MM)$-saturated model $(\NN, {\mathfrak X})$. But then $\NN \equiv \MM$, $\ssy(\NN) = \ssy(\MM)$, and $\NN$ is countable and recursively saturated, so $\NN \cong \MM$. We can then let $\NN = \MM$. Since ${\mathfrak X}$ is countable, we can let ${\mathfrak X} = \{U_0,U_1,U_2, \ldots\}$, and then  let $\MM^* = (\MM, U_0,U_1,U_2, \ldots)$. \qed

\bigskip

Having Proposition~4.2, we see that the following theorem implies Theorem~4.

 \bigskip 
  
  {\sc Theorem 4.3}: {\em If 
   $\MM^*$ is countable and  recursively supersaturated, then  there is $\NN^* \succ \MM^*$ such that $\Ltr(\NN^* / \MM^*) \cong ({\mathbf N}_5, \nu_3)$.}
  
  \bigskip
  
  {\it Proof}. For $n  \in M$, let $\alpha_n^{\MM^*} : {\mathbf N}_5 \into \eq\big((n+2) \times M^{n+1}\big)$ 
   be the function obtained by interpreting Definition~2.1 within $\MM^*$. Then, $\alpha_n^{\MM^*}$ is an $\MM^*$-representation of ${\mathbf N} _5$. Let ${\mathcal C}$ be the set of those $\MM^*$-representations $\alpha$ such that for some nonstandard $n \in M^*$,  $\MM^*\models \alpha \cong \alpha_n$.  It  is consequence of Theorem~2.9 and the recursive supersaturation of $\MM^*$ that ${\mathcal C}$ is an ${\MM^*}$-correct set (see Definition~$1.6^*$) of representations of $({\mathcal N}_5, \nu_3)$. Since $\MM^*$   is countable, Theorem~${1.7^*(2)}$ can be applied, yielding $\NN^* \succ \MM^*$ such that 
   $\Ltr(\NN^* / \MM^*) \cong ({\mathbf N}_5, \nu_3)$.    \qed

  \bigskip
  
  Theorem~4 now follows from Theorem~$1.7^*$.

% =======================================
% ------------------------------ REFERENCES-----------------
% ---------------------------------------------------------------------

\bibliographystyle{plain}

\end{document}